\documentclass[a4,11pt]{article}
\renewcommand{\baselinestretch}{1.3}
\newtheorem{propo}{Proposition}[section]

\newtheorem{theorem}{Theorem}

\newtheorem{lemma}{Lemma}[section]

\newtheorem{assumption}{Assumption}

\usepackage{xcolor}
\def\qed{ \ \vrule width.2cm height.2cm depth0cm\smallskip}
\def \wd {\wedge}

\def \rw {\rightarrow}

\newcommand{\eps}{\varepsilon}

\newcommand{\brm}{\begin{rem}}
\newcommand{\ermq}{\end{rem}}

\newcommand{\ba}{\begin{array}}
\newcommand{\ea}{\end{array}}
\newcommand{\be}{\begin{equation}}
\newcommand{\ee}{\end{equation}}
\newcommand{\bea}{\begin{eqnarray}}
\newcommand{\eea}{\end{eqnarray}}
\newcommand{\beaa}{\begin{eqnarray*}}
\newcommand{\eeaa}{\end{eqnarray*}}

\def \T {\Theta}
\def \R{I\!\!R}
\def \N{{\bf N}}
\def \bp {{\bf P}}

\def\g{\gamma}

\def\t{\tau}

\def\o{\omega}

\def\O{\Omega}
\def\cF{{\cal F}}

\def\cJ{{\cal J}}

\def\cP{{\cal P}}

\def\cT{{\cal T}}

\def\no{\noindent}

\def\ms{\medskip}

\def\q{\quad}
\def\qq{\qquad}

\def\bF{{\bf F}}

\def\qed{ \hfill \vrule width.25cm height.25cm depth0cm\smallskip}

\newcommand{\basa}{\begin{assumption}}
\newcommand{\easa}{\end{assumption}}

\newcommand{\bas}{\begin{assum}}
\newcommand{\eas}{\end{assum}}

\def\limsup{\mathop{\overline{\rm lim}}}

\def\esssup{\mathop{\rm esssup}}

\def\dis{\displaystyle}

\def\bF{{\bf F}}

\def \E {{\bf E}}

\newtheorem{thm}{Theorem}[section]
\newtheorem{lem}[thm]{Lemma}
\newtheorem{cor}[thm]{Corollary}

\newtheorem{rem}[thm]{Remark}

\newtheorem{defn}[thm]{Definition}
\newtheorem{assum}[thm]{Assumption}

\title{$\varepsilon$-Nash Equilibria of a Multi-player Nonzero-sum Dynkin Game in Discrete Time
}\author{Said
Hamad\`ene\thanks{Universit\'e du Maine, LMM, Avenue Olivier
Messiaen, 72085 Le Mans, Cedex 9, France. e-mail:
hamadene@univ-lemans.fr} , \, Mohammed
Hassani\thanks{Universit\'e Cadi Ayyad, Facult\'e poly-disciplinaire
de Safi, D\'epartement de Math\'ematiques et Informatique. B.P. 4162
Safi Maroc. e-mail : medhassani@ucam.ac.ma.
  This work has been carried out while the second author was visiting Universit\'e du
Maine, Le Mans (Fr.).} \,\,and 
Marie-Am\'elie Morlais$^*$
}

\begin{document}
\date{\today}
\maketitle
\begin{abstract}
We study the infinite horizon discrete time N-player nonzero-sum Dynkin game ($N\geq 2$)
with stopping times as strategies (or pure strategies). We prove existence of an $\eps$-Nash equilibrium point for the game by presenting a constructive algorithm. One of the main features is that the payoffs of the players depend on the set of players that stop at the termination stage which is the minimal stage in which at least one player stops. The existence result is extended to the case of a nonzero-sum game with finite horizon. Finally, the algorithm is illustrated by two explicit examples in the specific case of finite horizon. \end{abstract}
{\bf AMS Classification subjects}: 91A15 ; 91A10 ; 91A30 ; 60G40
91A60.
\medskip

\no {$\bf Keywords$}: Nonzero-sum Game ; Dynkin game ; Snell
envelope ; Stopping time ; Nash equilibrium point ; Pure strategy.

\section{Introduction}
The following zero-sum game on stopping times was introduced by 
Dynkin (1969).
Two players (or decision makers) $\pi_i$,
$i=1,2$, observe a bivariate sequence of adapted random variables
$\{(x_n,y_n), n\geq 0\}$. The first (resp. second) player chooses a stopping time which is denoted by
$\t_1$ (resp. $\t_2$) such that for any $n\geq 0$, $\{\t_1=n\}\subset \{x_n\geq 0\}$ (resp. $\{\t_2=n\}\subset \{x_n< 0\}$). At $\t_1\wedge \t_2$, if it is finite, $\pi_2$ pays $\pi_1$ an amount which equals to $y_\t$ and the game terminates.
If the game never terminates, $\pi_2$ does not pay anything. The objective of $\pi_1$ (resp. $\pi_2$) is to maximize
(resp. minimize) the following expected payoff
$$
\g(\t_1,\t_2)=\E [y_\t].$$  In Dynkin (1969), 
the author proved that if
$\sup_{n\geq 0}|y_n|$ is integrable, the game has a value, i.e.,
$$
\sup_{\tau_1}\inf_{\t_2}\g(\t_1,\t_2)=\inf_{\t_2}\sup_{\tau_1}\g(\t_1,\t_2).
$$ Moreover he also characterized $\varepsilon$-optimal stopping times. Since this seminal work, the discrete time zero-sum game has been widely discussed in several settings and works (Heller 2012; Kiefer 1971; Neveu 1975; Rosenberg et al 2001; Yasuda 1985). Such a list is far from exhaustive. 

Comparatively nonzero-sum Dynkin games have been less discussed, even if there are also several works on this subject (see e.g. Ferestein 2005; Nowak and Szajowski 1999; Hamad\`ene and Hassani 2013; Morimoto 1986; Neumann et al 2002; Ohtsubo 1987 and 1991; Shmaya et al 2003 and 2004).
 However, almost all of those works either deal only
with the case of two players and/or suppose some special structure of
the payoffs, or, the strategies of the players are of randomized
type. Therefore the main objective of our work is to study the discrete time nonzero-sum Dynkin
game when:\\
(i) there are more than two players and the strategies of players
are pure or
stopping times ;\\
(ii) the reward of each player, which is a stochastic process,
depends also on the set of players
which make the decision to terminate the game ;\\ 
(iii) the payoff processes are not supposed to satisfy a specific structure condition like being supermartingales or other structures (see Mamer 1987; Morimoto 1986). 

\ms In this paper, the problem we deal with is of nonzero-sum type in discrete time and it can be briefly described as follows. 

Let us consider $N$ players which are denoted
$\pi_1,...,\pi_N$ ($N\geq 2$). Let also assume that 
for $i=1,\dots, N$, $\pi_i$ chooses the stopping time $\t_i$ in order to stop or exit from the game which terminates at $R:=\min\{\t_j, j\in\{1,\cdots,N\}\}=\t_1\wedge ...\wedge \t_N$. The
corresponding payoff for $\pi_i$ is given by \be
\label{eqpayoff}\ba{lll} \dis J_{i}(\t_1,\cdots, \t_N)(\omega) :=
X^{{i},I(\omega)}_{R(\omega)}(\omega) \ea \ee where:

(i) $I(\omega)=\{j\in\{1,\cdots,N\}, 
\t_j(\omega)=R(\omega)\}$ is the coalition of players which make the decision to stop the game at $R(\omega)$ ;

(ii) $X^{i,I}$ is the stochastic payoff process for $\pi_i$ which depends on $I$. \ms

\noindent Our main objective is to find an $\eps$-Nash
equilibrium point (hereafter {\bf NEP} for short) for the game,
i.e., an $N$-tuple of stopping times $(\t_1^*,...,\t_N^*)$ such that
for any $i=1,...,N$, $$\eps+{\bf E}[J_i(\t_1^*,\cdots, \t_N^*)]\geq {\bf E}[J_i(\t_1^*,...,\t_{i-1}^*,
\t,\t_{i+1}^*,...,\t_N^*)],\,\,\mbox{for any stopping time }\t.$$
We mention here that the case of $\eps=0$ has been already considered in Hamad\`ene and Hassani (2013).

 \ms

In this paper, we show that the nonzero-sum discrete time game described above has an $\eps$-Nash
equilibrium point in pure strategies. It is a continuation of the work 
 of Hamad\`ene and Hassani (2013) 
where they have shown that the game has an $0$-{\bf NEP} if the payoff processes satisfy some specific property at infinity (see (\ref{onep}) below). Since we do not suppose this property, we cannot expect an $0$-{\bf NEP} for the game but only $\eps$-Nash equilibria. \ms

This paper is organized as follows: in Section 2, we set accurately
the problem, recall the Snell envelope notion and provide a result
(Theorem \ref{thmsnell}) which is in a way the streamline in the
construction of the $\eps$-{\bf NEP} for the discrete time nonzero-sum
Dynkin game we consider. We also discuss the relevance of the main assumption on the payoff processes (referred later as Assumption (A)) through two examples. The approximating scheme and its main properties are introduced in Section 3. In Section 4, we show that the
limit of the approximating scheme provides an $\eps$-{\bf NEP} for the game which is the main result of the paper. We also provide an extended result to the case of nonzero-sum games with finite horizon. Finally, Section 5 is devoted to the analysis of two examples of games with finite horizon: the first with deterministic rewards and the second one with stochastic rewards. For both examples, the constructive algorithm is carried out. We  shall provide explicit $\eps$-{\bf NEP} and discuss some of their properties such as the existence of several $\eps$-Nash equilibria.

\section{Setting of the problem and hypotheses}
\setcounter{equation}{0} Throughout this paper, $\N$ denotes the set of integers and $(\Omega
,\mathcal{F},\bp)$ is a
 fixed probability space on which is defined a filtration
 ${\bF}:=(\cF_{t})_{t\in\N}$. For any stopping time $\theta $, let us denote by
 \medskip

(i) $\cT_\theta$ the set of all $\N$-valued $\bF$-stopping times $\tau$ such that $\tau \geq\theta$ and $\cT_\theta^N=\underbrace{\cT_\theta\times ...\times \cT_\theta}_{N \mbox{ times }}$ ;

(ii) $\E_\theta[.]$ the conditional expectation with respect to ($w.r.t.$ for short) $\cF_\theta$, i.e., ${\bf E}_\theta[X]:={\bf E}[X|{\cal F}_\theta]$,
for any integrable random variable $X$;

(iii) $\cJ:=\{1,...,N\}$ and $\cP := \{I\subseteq \cJ \mbox{ such
that } I\neq\emptyset \}$. \bigskip

Let $\bar{\N}:=\N\cup \{\infty\}$ and without loss of generality (or w.l.o.g in short), we assume that ${\cal F}_\infty:={\cal
F}=\bigvee_{t\geq 0}{\cal F}_t$. For $i\in\cJ $ and $I\in\cP$, let
$(X^{i,I}_t)_{t\in\bar{\N}}$ be an ${\bF}$-adapted and real valued process
such that $${\bf E}[\displaystyle\sup_{t\in\bar{\N}}
|X^{i,I}_t|]<\infty.$$ We moreover assume that they satisfy the
following hypotheses.
\begin{assum}
\label{assumption}${}$ $\bf {(A)}:$ For any $i,j=1,\dots,N$ and all $t\in \N$,
$$ X_t^{i, \{i,j\}}\leq X_t^{i,\{j\}},\,\,\bp-a.s.\qq\qed$$
\end{assum}

For $T_1, \cdots, T_N$ elements of $\cT_0$ and $i\in
\cJ$, we define $J_i(T_1,T_2,\cdots, T_N)(\omega)$, the payoff associated
with the player $i$, as follows: \be \label{J12} \ba{lll} \dis
J_i(T_1,T_2,\cdots, T_N) := \sum_{I\in\cP}X^{i,I}_{R}\,1_{\cap_{j\in
I}\{T_j=R \}\cap\cap_{j\in I^c}\{T_j> R \}},\,\bp -a.s. \ea \ee where:

(i) $R:=\min\{T_j, j\in\cJ\}=T_1\wedge ... \wedge T_N$ ; 

(ii) by convention, we assume that $\cap_{i\in\emptyset}A_i=\Omega$.
\medskip

\no Let us emphasize that, for fixed $\omega$, if ${I}(\omega):= \{j\in\cJ,  \,
T_j(\omega)=R(\o) \}=I_0$ then
$$\dis J_i(T_1,T_2,\cdots, T_N)(\omega) =X^{i,{I_0}}_{R}(\omega) .$$
Note that if $R(\o)=\infty$ then obviously $I(\o)=\cJ$.
\ms
\begin{rem} If $I_0\neq\cJ$ then $X^{i,I_0}_\infty$ does not play any role, therefore and w.l.o.g. we can assume that  $X^{i,I_0}_\infty=0$.
\end{rem}

We next precise the notion of $\varepsilon$-equilibrium we deal with.
\begin{defn}
\label{equidefn} Let $\varepsilon\geq 0$. An $N$-tuple of stopping times $(T_1^*,T_2^*,\cdots, T_N^*)\in
{\cT_0}^N$ is a called an $\varepsilon$-{\bf NEP} point for the nonzero-sum Dynkin game if for all $i=1,\cdots,N$ we have
 \be
\label{equilibrium} {\bf
E}[J_i(T_1^*,\cdots,T_{i-1}^*,T,T_{i+1}^*,\cdots, T_N^*)]\le {\bf
E}[J_i(T_1^*,\cdots,T_{i-1}^*,T_i^*,T_{i+1}^*,\cdots, T_N^*)]+\varepsilon,
\,\forall \; T\in {\cT_0}. \ee
\end{defn}

\begin{rem} (i) If $\varepsilon =0$, this definition means that $(T_i^*)_{i=1,N}$ is a standard {\bf NEP} for the game. Otherwise, i.e., if $\epsilon >0$, it means that for any $i=1,\dots, N$, $(T_i^*)_{i=1,N}$ verifies:
$$
\left| \sup_{T\in {\cT_0}}{\bf
E}[J_i(T_1^*,\cdots,T_{i-1}^*,T,T_{i+1}^*,\cdots, T_N^*)]-{\bf
E}[J_i((T_i^*)_{i=1,N})]  \right|  \leq \eps. $$ 
\ms

\noindent (ii) In order to show that the game has an $\eps$-{\bf NEP}, we need Assumption (A) to be fulfilled. However we do not know how to get rid of it since, when it is not satisfied, the game may or may not have an $\eps$-{\bf NEP}. This can be seen through the two following examples. 

Actually assume that for any $n\geq 0$, $\cF_n=\{\Omega, \emptyset\}$. Then $\cT_0$ is reduced to constant stopping times. Next for $n\in \bar \N$, let us set:
$$ 
X^{1,\{1\}}_n=0, X^{1,\{2\}}_n=0, X^{1,\{1,2\}}_n=1 \mbox{ and }X^{2,\{2\}}_n=0, X^{2,\{1\}}_n=0, X^{2,\{1,2\}}_n=-1.$$ 
Then the assumption (A) is not satisfied since 
$$ \forall n\in \bar \N, \;\; \;
X^{1,\{1,2\}}_n=1>0=X^{1,\{2\}}_n.$$  
 On the other hand we have, $$ J_1(t_1,t_2)=-J_2(t_1,t_2)=1_{(t_1=t_2)} .$$  Therefore one can easily check that for $\eps$ in $(0,1)$ this nonzero-sum Dynkin game does not have an $\varepsilon$-$ {\bf NEP}.$ 
 
 \ms 
Let us now skip to the following second example. For $n\in \bar \N$, let us define $$ 
X^{1,\{1\}}_n=0, X^{1,\{2\}}_n=0, X^{1,\{1,2\}}_n=1 
\mbox{ and }X^{2,\{2\}}_n=0, X^{2,\{1\}}_n=0, X^{2,\{1,2\}}_n=0.$$ 
Thus and once more, (A) is not satisfied since 
$$
\forall n\in \bar \N, \; \; \; X^{1,\{1,2\}}_n=1>0=X^{1,\{2\}}_n. $$
On the other hand $$ J_1(t_1,t_2)=1_{(t_1=t_2)} \hbox{ and } J_2(t_1,t_2)=0.$$
Then, for any $\varepsilon\geq 0$ and $\underbar t$ arbitrarily fixed in $\N$, $(\underbar t,\underbar t)$ is an $\eps$-$ {\bf NEP}$, which means that $(A)$ is not a necessary condition. 
\ms

\noindent (iii) Under Assumption (A), if moreover the processes $X^{i,I}$ verify  
\be\label{onep}\displaystyle \limsup_{t\in
\N}X_t^{i,\{i\}}
:=
\inf_{t\in\N} \sup_{t\leq n <\infty}X_n^{i,\{i\}}\leq
X_\infty^{i,\cJ},\, \;\,\bp-a.s.\ee
 then it is proved in Hamad\`ene-Hassani (2013) that the game has an $0$-{\bf NEP}.
 \qed
\end{rem}

To tackle the game problem we consider, we mainly use the notion of Snell envelope of processes which we introduce briefly below. For more
details on this subject, one can refer either to Dellacherie-Meyer (1980), pp.431 or to El Karoui (1980), pp.140. 
For sake of completeness we give the
following result related to existence of an $\eps$-optimal stopping time as
we do not find a reference where it is given in the form we need later.
\begin{theorem}\label{snev}  \label{thmsnell} 
 Let $\;U=(U_t)_{t\in \bar{\N}}$ be an
$\bF$-adapted $\R$-valued  process such that ${\bf E}[\sup_{t\in
\bar{\N}} |U_t|]<\infty$. For any $\bF$-stopping time $\theta$
let us define:  \be \label{sun} Z(\theta)=\esssup_{\tau \in {\cal
T}_{\theta}}{\bf E}[U_\tau|\cF_\theta]\, (\mbox{and then
}Z(\infty)=U_\infty).\ee For $n\in \bar{\N}$, let us set
$Z_n:=Z(n)$. Then

(i) $Z(\theta)=Z_\theta$, $\forall \theta \in \cT_0$ ;

(ii) $(Z_n)_{n\geq 0}$ is an $\bF$-supermartingale which satisfies:
$$
\forall t\ge 0,\,\,Z_t=U_t\vee E_t[Z_{t+1}], \; \bp-a.s.
$$
Moreover and for any $\varepsilon >0$, the stopping
time $ \tau^*$ defined by
$$\tau^*=\inf\{s\geq0,\q Z_s\leq U_s+\varepsilon\},$$ is $\varepsilon$-optimal, $i.e.$,
\begin{equation}\label{sdeux} 
\sup_{\tau \in \cT_0}
{\bf E}[U_\tau]\leq {\bf E}[U_{\tau^*}]+\varepsilon.\end{equation}
Finally, $\lim_{t\rw \infty}Z_t= U_\infty
 \hbox{ on the set } \big( \tau^*=\infty\big).$
\end{theorem}
{\bf Proof}: First note that for any stopping time $\theta$, the random variable $Z(\theta)$ is defined since ${\bf E}[\sup_{t\in
\bar{\N}} |U_t|]<\infty$. Next the first property (i) follows from the fact that, for all stopping times $\theta$ and $\lambda$, $$Z(\theta)=Z(\lambda) \hbox { on the random set } \{\theta=\lambda\},$$ and since $\bar{\N}$ is
a discrete set. Let us focus on (ii). For any $t\in \N$ we
have,
$${\bf E}_t[Z_{t+1}]=\esssup_{\tau \in {\cal T}_{t+1}}{\bf
E}[U_\tau|\cF_t]\leq Z_t, $$ which implies that 
$(Z_n)_{n\geq 0}$ is an $\bF$-supermartingale. On the other hand and by definition
$$Z_t\geq U_t \vee
E_t[Z_{t+1}].$$  For any $ \tau\in{\cal T}_{t}$, it also holds $$  {\bf E}[U_\tau|\cF_t] =  U_t 1_{(\tau=t)} + {\bf
E}[U_{\tau\vee (t+1)}|\cF_t]1_{(\tau\geq t+1)}\leq U_t
1_{(\tau=t)} + {\bf E}[Z_{t+1}|\cF_t]1_{(\tau\geq t+1)} \leq
U_t \vee E_t[ Z_{t+1}],$$ which implies that $$  \forall \; t \in \N, \q Z_t = U_t \vee
E_t[ Z_{t+1}].
$$
Since $(\tau^*>  t)\subset  (Z_{t}>U_{t}+\eps)$, we claim $$\forall \;t \in \N,  \q ({\bf E}[Z_{(t+1)\wedge\tau^*}|\cF_t]-Z_{ t\wedge\tau^*})= \big({\bf
E}[Z_{t+1}|\cF_t] - Z_t\big)1_{(\tau^*>t)} = 0.$$
 Thus, it gives
\begin{equation}\label{pouf}{\bf E}[Z_{(t+1)\wedge\tau^*}]={\bf
E}[Z_{t\wedge\tau^*}]= {\bf
E} [Z_0].
\end{equation}
The supermartingale 
$(Z_{t\wedge\tau^*})_{t\in \N}$ is actually a martingale. Besides we have
\begin{equation} \label{unifint}\forall n\geq 0,\,|Z_n|\leq {\bf E}[\sup_{k\in \bar \N}|U_k||\cF_n]. \end{equation} Henceforth: 
  
(a) the supermartingale 
$(Z_{t})_{t\in \N}$ is $\bp$-a.s. convergent and uniformly integrable and it converges in $L^1(d\bp)$ ;  

(b) the martingale $(Z_{t\wedge\tau^*})_{t\in \N}$ is uniformly integrable and then converges in $L^1(d\bp)$ to $Z_{\t^*}$; 

(c) the random variable $Z_{\t^*}$ is integrable and by Fatou's Lemma we have $ \E[Z_{\t^*}] \leq \E[Z_0]$.
\ms

\noindent But for any $t\geq
s$ and $\tau \in {\cal T}_{t}$ we have $${\bf E}[U_\tau|\cF_t]\leq
{\bf E}[\sup_{s\leq n\leq\infty} U_n|\cF_t],$$ and then $$Z_t\leq {\bf
E}[\sup_{s\leq n\leq\infty} U_n|\cF_t].$$ 
Therefore taking the limit in $t$ and using a result by Neveu (1975) (see  Proposition II.2.11 pp.29) we obtain
$$
\lim_{t\rw \infty}Z_t\leq \sup_{s\leq n\leq\infty} U_n, $$
from which we get, by taking the infimum w.r.t $s$, the following inequality
\be \label{eqint1}
\lim_{t\rw \infty}Z_t\leq \limsup_{t\rw \infty} U_t\vee U_\infty.\ee
\noindent Next by taking the limit in $t$ in (\ref{pouf}) and taking into account (a), (b) and (c) above yields  \be \label{resopti}\begin{array}{ll}0\leq 
{\bf
E}[Z_0] -{\bf E}[Z_{\tau^*}] &= \lim_{t\rw \infty}{\bf E} [Z_{ t\wedge\tau^* }-Z_{\tau^*}]
\\ &= {\bf E}[\big(\lim_{t\rw \infty}Z_{t}-U_{\infty}\big)1_{(\tau^*=\infty)} ].\end{array}\ee
 But on $(\tau^*=\infty)$ we have  $$ \lim_{t\rw \infty}Z_t\geq \limsup_{t\rw \infty} U_t
 +\varepsilon.$$
Then and using (\ref{eqint1}) we have $$ \limsup_{t\rw \infty} U_t+\varepsilon \leq \lim_{t\rw \infty}Z_t\leq U_\infty.$$  (\ref{resopti}) finally implies that $$ {\bf
E}[Z_0]={\bf E} [Z_{\tau^*}]\q \hbox{ and }\q \bp-a.s., \,\lim_{t\rw \infty}1_{\{ \tau^*=\infty\}}Z_t= U_\infty 1_{\{ \tau^*=\infty\}}.$$  
  Thus
$$\sup_{\tau\geq0}{\bf E}[ U_\tau]= {\bf
E}[Z_0] ={\bf E}[Z_{\tau^*}]\leq {\bf
E}[U_{\tau^*}]+\varepsilon, $$ which means that $\tau^*$ is $\varepsilon$-optimal.  
\qed

\begin{rem} (i) $(Z_n)_{n\in \bar \N}$ is actually the smallest ${\bF}$-supermartingale which is greater than the payoff process $U$. 

\noindent (ii) If the condition $\limsup_{n\rw \infty}U_n\leq U_\infty$, is not satisfied then an 0-optimal stopping time may not
exist. To illustrate this claim, let us consider the process $(U_n)_{n\in \bar \N}$ defined as follows  \\$U_n = 1 -\frac{1}{n+1} \mbox{ for any } n \mbox{ in } \N \mbox{ and }U_\infty=0.$
Due to the first point (i) and using (\ref{eqint1})  $$ U_n = 1 - \frac{1}{n} \le Z_n \le 1 = \limsup U_n, $$ and thus, the Snell envelope $(Z_n)_{n\in \bar \N}$
of $U$ is defined as follows $$ Z_n = 1, \; \mbox{ if }n \in \N \;\; \mbox{ and }   \; \; Z_\infty = 0.$$
Now if $\t$ is a stopping time then
$${\bf
E}[U_\t] = {\bf
E}[U_\t1_{(\t<\infty)} ] = \bp[\t<\infty]-{\bf
E}[\frac{1}{\t+1}]<1=\sup_{\t \in {\cal T}_0}{\bf E}[U_\t]={\bf
E}[Z_0].$$
Thus an 0-optimal stopping time does not exist for the optimal stopping problem with payoff $U$. However and for any $\eps >0$, if $n_\eps$ is such that $\frac{1}{n_\eps+1}<\eps$ then $n_\eps$ is an $\eps$-optimal stopping time. 
\end{rem}

\section{The approximating scheme and its properties}

Let us introduce sequences of stopping times which, as it will be shown later, converge to an $\varepsilon$-{\bf NEP} of the game. Hereafter, $\varepsilon>0$ is fixed and we define by induction a sequence of $\bF$-stopping times $
(\tau_n)_{n\geq1}$ in circular way since there is a move from one player to the next one until all the objects are defined for all players. Then, the procedure starts again with the first player. More precisely and for $n\geq 1$, let $(i_n,q_n)$ be the unique pair of integers such that
$n=Nq_n+i_n$ with $i_n\in\{1,2,\cdots,N\},$ and let us set: \ms
\setcounter{equation}{0} 
\\
\no (i) $\tau_1=\cdots\tau_N=\infty,$ and\\
(ii) for $n\geq N+1$,  we put   $$ 
\begin{array}{ll} 
(a)\q \theta_n= \min\{\tau_{n-1},\tau_{n-2},\cdots,
\tau_{n-(N-1)}\}\,;\\ \\
(b)\q I_n:=\{i_l\in\cJ : n-N+1\leq l \leq n-1 \hbox{ and }
\tau_l=\theta_n\}\, ;\\ \\
 (c) \q  \forall t\in\bar{\N}, \,U^n_t= X^{i_n,\{i_n\}}_{t}1_{\{t<\theta_n \}} +
 Y^n 1_{\{t\geq\theta_n \}} \hbox{ with } 
   \\\\ \qq \qq
    Y^n=\big(X^{i_n,I_n\cup\{i_n\}}_{\theta_n}\vee X^{i_n,I_n }_{\theta_n}\big) 1_{\{\theta_n<\infty\}} + X^{i_n,\cJ}_{\infty} 1_{\{\theta_n  =\infty\}}\,;
    \\ \\
(d) \q  \forall t\in\bar{\N}, \, W^n_t=\esssup_{\nu \in {\cal T}_{t}}{\bf E}[U^n_\nu|\cF_t]\,;\\ \\
(e) \q  \mu_n=\min\{s\in\bar{\N}, W^n_s\leq U^n_s+\varepsilon\}\, ;\\ \\
(f) \q \tau_n=
(\mu_n\wedge\tau_{n-N})1_{\{\mu_n\wedge\tau_{n-N}<\theta_n\}}+\tau_{n-N}1_{\{\mu_n\wedge\tau_{n-N}\geq\theta_n\}}.
\end{array}
$$
A few properties are collected below in the following remark.
\begin{rem} \label{rem1} 
For any $n\geq N+1$,\\
\no (i) $i_n$ does not belong to $I_n$ and for every
$I\in\cP$ such that $i_n\notin I$ we have
$$\{I_n=I\}:=\{\o\in \O, I_n(\o)=I\}= \cap_{j\in I} (\tau_{k_j}=\theta_n)\cap\cap_{j\in I^c\backslash\{i_n\}} (\tau_{k_j}>\theta_n)\in \cF_{\theta_n},$$
where, for $j\neq i_n$, $k_j $ is the unique integer such that $k_j\in\{n-N+1,\cdots,
n-1\}$ and $i_{k_j}=j.$

 \no (ii) $W^n$ is a supermartingale that satisfies for all $t\geq\theta_n$
$$ 
W^n_{t} =U^n_{t}=Y^n. $$ 
Moreover the process $(W^n_{t\wedge \mu_n})_{t\geq 0}$ is an $\bF$-martingale. \\
(iii) The following inequalities are satisfied: \be\label{ineq2}  \mu_n \leq \theta_n,\,\,
\tau_n\leq\tau_{n-N}\mbox{ and } \theta_{n}\leq\theta_{n-N}. \ee
(iv) By Theorem \ref{snev}-(ii), the stopping time $\mu_n$ is $\varepsilon$-optimal, i.e.,
$$ \forall\tau \in \cT_0, \,\, {\bf E}[U^n_\tau]\leq {\bf E}[U^n_{\mu_n}]+\varepsilon. \,$$
(v) Let $n_0$ be fixed. Since the induction is of circular type, then player $i_{n_0}$ knows that the game will be terminated at $\theta_{n_0}$ and $U^{n_0}$ is her payoff. She then chooses the time $\t_{n_0}$ to stop the game accordingly. \qed
\end{rem}

\noindent First we are going to simplify the expression of $\t_n$. 
\begin{propo}\label{simplif}
For any $n\geq 1$, $ \mathbf{P}( \{ \mu_{n+N}\leq \t_n \}) = 0$.
\end{propo}
{\bf Proof}: Suppose on the contrary that there exists $m\geq 1$ such that
$\mathbf{P}[\tau_m < \mu_{m+N}]>0$. Let us set $n=\min\{m\geq1 \;
\hbox{s.t.}\; P[\tau_m < \mu_{m+N}]>0\}$. 
Since $\tau_1= \cdots =\tau_N =\infty$, then necessarily
$n\geq N+1$. On the set $\Theta:=\{\tau_n < \mu_{n+N}\}$, such that $\mathbf{P}(\Theta) >0$ by definition we
have
\begin{equation}\label{eq4}\tau_n < \theta_{n+N}:= \tau_{n+N-1}\wedge
\tau_{n+N-2}\wedge\cdots \tau_{n+1},\end{equation}since
$\mu_{n+N}\leq \theta_{n+N}$ (see Remark \ref{rem1}-(ii)). Thus the minimality 
$n$ implies that for all $j$ such that 
$j\in\{1,2,\cdots,n-1\} $,
$\mu_{j+N}\leq \tau_{j}$ and by definition of $\tau_{j+N}$
$$\tau_{j+N}=\mu_{j+N} 1_{\{\mu_{j+N}
<\theta_{j+N}\}}+\tau_{j}
1_{\{\mu_{j+N}=\theta_{j+N}\}},$$
since $ \mu_{j+n} \wedge \tau_j  = \mu_{j+N}.$ From
(\ref{eq4}) and the definition of $\theta_{n+N-1}$ we deduce that
$\theta_{n+N-1}=\tau_n $ on $\Theta$. It follows that
\begin{equation}\label{newdef}
\tau_{n+N-1}=\mu_{n+N-1} 1_{\{\mu_{n+N-1}
<\t_n\}}+\tau_{n-1}1_{\{\mu_{n+N-1}=\t_n\}}.
\end{equation}
Therefore, once more on $\Theta$, we claim that
\be\label{eq5} \tau_n <\tau_{n+N-1} =\tau_{n-1}. \ee The strict
inequality in (\ref{eq5}) stems from (\ref{eq4}). Noting that $\mu_{n+N-1}\le \theta_{n+N-1}=\t_n$ \\
and $\tau_n <\tau_{n+N-1}$ on $\Theta$, we obtain $ \mu_{n+N-1} <\tau_{n+N-1}$. Combined with (\ref{newdef}), the equality in (\ref{eq5}) holds true.\\
Let us now justify the following property on the set $\Theta$ 
\begin{equation}\label{stationarity}
 \forall \; j \in \{1,\dots, N-1\}  \; \quad \tau_{n-j} = \tau_{n+N-j}.
\end{equation}
 Since the claim already holds for $j=1$, we prove it for $j=2$. 
By definition of $\Theta_{n+N-2}$, one has 
 $$\theta_{n+N-2}= \tau_{n+N-3}\wedge \tau_{n+N-4}\wedge\cdots
\tau_{n+1}\wd \t_n\wd \t_{n-1} =\t_n.$$ Indeed, using first (\ref{eq4}), we obtain $\t_n<\t_{n+k}$,  for
any $k$ in $1,...,N-1$ and using (\ref{eq5}), we claim that 
$\t_n<\t_{n-1}$. Using once again the minimality of $n$ and the definition of
$\tau_{n+N-2}$, we obtain 
$$
\tau_{n+N-2}=\mu_{n+N-2}
1_{\{\mu_{n+N-2}<\t_n\}}+\tau_{n-2}1_{\{\mu_{n+N-2}=\t_n\}}.
$$
Thanks to (\ref{eq4}), $\t_n < \tau_{n+N-2}$ and thus  
$\t_n < \tau_{n+N-2}=\tau_{n-2}$. Assuming otherwise that\\ $\tau_{n+N-2}=\mu_{n+N-2}$, it yields $\tau_{n+N-2} <\t_n $ which is a contradiction on $\Theta$ and gives us the desired result for $j=2$.
 Repeating the same arguments as many times as necessary, we obtain the claim stated in (\ref{stationarity}).\\
Therefore, due to property (\ref{stationarity}) and on the set $\Theta$, it holds 
$$
\t_n<\theta_{n+N}=\theta_n \hbox{ and } I_{n+N}=I_n.$$ Using both the minimality of $n$ and the definition of $\t_n$ we obtain 
$$\tau_n=\mu_{n} 1_{\{\mu_{n}
<\theta_{n}\}}+\tau_{n-N}1_{\{\mu_{n}=\theta_{n}\}}=\mu_{n},
$$ since $\mu_n\leq \t_{n-N}$ and $\t_n<\theta_n$. Henceforth on $\Theta$, 
we have $U^n=U^{n+N}$ since
$\theta_{n+N}=\theta_n$, $i_{n+N}=i_n$ and $ I_{n+N}=I_n$. By definition, we obtain
$$\begin{array}{lll} 1_\T W^{n+N}_{\mu_n}=1_\T
W^{n+N}_{\t_n}&=&1_\T\esssup_{\nu \in {\cal T}_{\t_n}}{\bf E}[
U^{n+N}_\nu|\cF_{\t_n}]=\esssup_{\nu \in {\cal T}_{\t_n}}{\bf
E}[1_\T U^{n+N}_\nu|\cF_{\t_n}]
\\ & =& \esssup_{\nu \in {\cal T}_{\t_n}}{\bf E}[1_\T U^{n}_\nu|\cF_{\t_n}]
\\ &=& 1_\T W^{n}_{\t_n}=1_\T W^{n}_{\mu_n}\leq 1_\T (U^{n}_{\mu_n}+\varepsilon)=1_\T (U^{n+N}_{\mu_n}+\varepsilon)
\end{array}
$$ i.e., $1_\T W^{n+N}_{\mu_n}\leq1_\T (U^{n+N}_{\mu_n}+\varepsilon)$ and then
$\mu_{n+N}\leq \mu_n$ on $\T$. As on $\T$ we have
$\mu_n=\t_n<\mu_{n+N}$, this is contradictory with the previous
inequality. Henceforth $\bp[\T]=0$ and for any $m\geq 1$ we have
$\mu_{m+N}\leq \t_m$, $\bp$-a.s., which completes the proof. \qed \ms

\noindent As a by-product, we obtain the following simplified expression of $\t_n$.

\begin{cor}\label{cor1} For any $n\geq N+1$,

(i) $ \tau_n=
\mu_n1_{\{\mu_n<\theta_n\}}+\tau_{n-N}1_{\{\mu_n=\theta_n\}}  \le \tau_{n-N}$;

(ii)
$\mu_n=\tau_n\wedge\theta_n=\tau_n\wedge\tau_{n-1}\wedge\cdots\tau_{n-N+1}\leq \mu_{n-N}
.$
\end{cor}
{\bf Proof}: Using both Proposition \ref{simplif} 
and the definition of $\t_n$, we obtain (i). As for (ii), 
for any $n\geq N+1$, we have
$$\ba{ll} \tau_n\wedge\theta_n&=\tau_n 1_{\{\tau_n< \theta_n\}}+
\theta_n 1_{\{\tau_n\geq \theta_n\}}.\ea$$But $\tau_n
1_{\{\tau_n< \theta_n\}}=\mu_n 1_{\{\mu_n< \theta_n\}} $ and
on $[\tau_n\geq \theta_n]$ we have $\theta_n=\mu_n$. Therefore
$\theta_n 1_{\{\tau_n\geq \theta_n\}}=\mu_n 1_{\{\t_n\geq
\theta_n\}}=\mu_n 1_{\{\t_n\geq \mu_n\}}$. Gathering now those
equalities yields $\mu_n=\tau_n\wedge\theta_n$. Finally the second
equality is just the definition of $\theta_n$. \qed
\ms

\noindent We state below some properties of the sequences $(\t_n)_n$, $(\theta_n)_n$, $(\mu_n)_n$, which we need later.

\begin{propo}\label{prop32}
For any $m\geq N+1$, $${\bp}[\{ \tau_m=\theta_m<\infty \}]=0.$$ 
\end{propo}
{\bf Proof}: Let $m\geq N+1$ and $\Omega'_m:=\{ \tau_m =\theta_m<\infty\}$. On the set $\O'_m$, we claim   \be\label{recu}\forall j\in\{m-(N-1),\cdots, m\},
 \q\theta_{j}\leq\tau_m \hbox{ and }
\tau_{j}=\tau_{j-N}.\ee To begin with, we note that for $j=m$ and by definition of $\tau_m$ in (i), Corollary \ref{cor1}
$$
\tau_{m}=\tau_{m-N}.$$
Indeed and on the set $\Omega'$,  $\tau_m = \mu_m$ is contradictory with $\tau_m = \Theta_m$. Property (\ref{recu}) is now proved for $j=m$. 
 We note that necessarily $m \ge 2N+1$. Indeed, since $\tau_{m-N} = \tau_m =\theta_m < \infty$ and since, by construction, $\tau_1= \cdots \tau_N =+ \infty $ then necessarily $ m-N \ge N+1$. 
Next we proceed with a backward induction procedure by supposing   
\begin{equation} \label{prop_taul} \exists \; l \in\{m-(N-1),\cdots,m-1\}, \;\forall \; 
j\in\{l+1,\cdots,m\}, \quad \theta_{j}\leq\tau_m\;\; \textrm{and} \; \;
\tau_{j}=\tau_{j-N}.
\end{equation}  Fixing $l$ satisfying (\ref{prop_taul}), we have to prove both $\theta_l\leq \t_m$ and $\t_l=\t_{l-N}$. By definition  $$\begin{array}{ll}\theta_{l}
&=\tau_{l-1}\wedge\cdots \wedge \tau_{m-N } \wedge\cdots\wedge\tau_{l-N+1}\\{}&=
\tau_{l-1}\wedge\cdots \wedge \tau_{m}\wedge
 \cdots\wedge\tau_{l+1}.\end{array}$$
The second equality follows from the induction hypothesis since $\t_{m-N}=\t_m$ and therefore, one obtains $\theta_l\leq \t_m$. Next, if $\t_{l}<\t_{l-N}$ and by definition of $\t_l$ we have $\t_{ l}<\theta_l\leq \t_m=\theta_m$.
Using once more the definition of $\theta_m$, 
we have $\theta_m\leq\tau_l$ since $l\in \{m-(N-1), \cdots,m-1\}$, which is absurd. Therefore 
$\t_l=\t_{l-N}$ and the proof of the induction is stated.\\
Relying now on (\ref{recu}) we have $$\begin{array}{ll}\theta_{m-N}&=\tau_{m-N-1}\wedge\cdots \wedge \tau_{m-N-N+1 }  \\{}&= 
\tau_{m-1}\wedge\cdots \wedge \tau_{m-N+1}= \theta_m.\end{array}$$ Since $\t_m=\t_{m-M}$, $  \Omega'_m \subseteq \Omega'_{m-N}$. Proceeding with an induction procedure, one gets
$$ \forall \; q \in \mathbf{N}, \; m-qN \in \mathbf{N}, \;\; \; \Omega'_m \subseteq \Omega'_{m-qN}.$$  
Thus and for $q$ large enough, $m-qN < N+1 $. Therefore and on the set $\Omega'_{m-qN}$, $\tau_{m-qN}  < \infty$ which is contradictory with $\tau_1 = \cdots \tau_N =\infty.$
This yields $ \Omega'_m  \subset \Omega'_{m-qN} =  \emptyset, \; \bp-a.s..$ The main claim  stated in Proposition \ref{prop32} is established. \qed

\begin{lemma} \label{mu} For any 
$m\geq N+1$, $$ \bigg( \mu_m=\mu_{m+N}\bigg)\subset\bigg( \t_m=\t_{m+N}\bigg).$$ 
\end{lemma}
{\bf Proof}: Let $m\geq N+1$. On the set $ \bigg( \mu_m=\mu_{m+N}\bigg)$ and assuming $\t_m>\t_{m+N}$ then
$$\mu_m=\mu_{m+N}=\t_{m+N}<\theta_{m+N}\leq\theta_m.$$
But since $\mu_m=\t_m\wedge\theta_m$, we have $\t_{m+N}=\mu_m=\t_m$ which is absurd and completes the proof. \qed


\section{Existence of an $\eps$-Nash equilibrium point for the game}
For any $i$ in $ \{1,\cdots,N\}$, let us define
\be \label{explicitNEP}
T_i^*=\displaystyle\lim_{n\longrightarrow\infty}\;\tau_{Nn+i}\; \hbox{
and }\;
R^*_i=\displaystyle\lim_{n\longrightarrow\infty}\;\theta_{Nn+i}=\min\{T_j^*
\; ; j\neq i\}.
\ee Those limits exist since for any $n\geq N+1$, we
know that $\t_n\leq \t_{n-N}$ therefore the sequences of stopping
times $(\tau_{Nn+i})_{n\geq 0}$ are non-increasing for any fixed
$i$. On the other hand, as $N$ is finite, we also have
$$\forall \; i\in
\{1,\cdots,N\},\,\,R^*:=T_1^*\wedge\cdots\wedge T_N^*=R^*_i\wedge
T^*_i=\displaystyle\lim_{n\longrightarrow\infty}\;\mu_{Nn+i}
=\displaystyle \lim_{n\longrightarrow\infty}\;\mu_n=\min\{\mu_n ; n\in\N\}.$$    Next for $i\in
\{1,\cdots,N\}$, let us define $$I^*_i(\o):=\{j\in
\cJ\backslash\{i\} : T^*_j(\o)=R^*_i(\o)\}.$$

In what follows, we show that the $N$-tuple of stopping times
$(T_i^*)_{i=1,...,N}$ is an $\eps$-{\bf NEP} point for the N-players
nonzero-sum Dynkin game associated with $(J_i)_{i\in \cJ}$. The 
proof is obtained after several intermediary results which involve the stationary decreasing sequences of stopping times $(\t_{nN+i})_{n\geq 0}$ and
their limits. For clarity, we list below the main steps: 

(i) we first establish a link between the payoffs
$J_i(T_1^*, \dots,T_{i-1}^*,\theta,T_{i+1}^*,\dots,T_{N}^*)$ and \\$
\lim_{n\rw \infty}{\bf
E}[ U^{Nn+i}_{\theta\wedge\theta_{Nn+i}}].
$
The stationarity of the sequences plays an important role here.

(ii) By using the link between $U^n$ and its Snell envelope process $W^n$, which is commonly used in optimal stopping problems, we are
able to compare $J_i(T_1^*, \dots,T_{i-1}^*,\theta,T_{i+1}^*,\dots,T_{N}^*)$ and $J_i(T_1^*,\dots, T_{N}^*)$ for any given fixed stopping time $\theta$.

(iii) Relying on Assumption (A), it allows us to cancel some extra terms and to check that $(T_i^*)_{i=1,\dots, N}$ is actually an $\eps$-{\bf NEP} for the game. 

\begin{lem}\label{lem1}
Let $(\beta_n)_{n\geq1}$ be a decreasing sequence of stopping times that converges to $\beta$.
Then for any $i\in\{1,\cdots,N\}$ 
we have \be \label{ebauche1}\begin{array}{ll} 
 \dis{ \lim_{n\rw \infty} }{\bf
E}[ U^{Nn+N+i}_{\beta_n\wedge\theta_{Nn+N+i}}] = \\ \\
\q   {\bf
E}[J_i(T^*_1,T^*_2,\cdots,T^*_{i-1},\beta,T^*_{i+1},\cdots,
T^*_{N})]
+\\\\
 \qq {\bf
E} [\bigg(X^{i,I^*_{i}}_{R^*_{i}}-X^{i,I^*_{i}\cup\{i\}}_{R^*_{i}} \bigg)^+ 1_{\{R^*_{i}=\beta<\infty\}}
+\bigg(X^{i,I^*_{i}}_{R^*_{i}}-X^{i,I^*_{i}\cup\{i\}}_{R^*_{i}} \bigg)^- 1_{\{R^*_{i}<\beta\}}].
\end{array}\ee
\end{lem}

\noindent {\bf Proof}:
First and since $\beta_n$ is $\mathbf{N}$-valued we obtain that the sequence 
$(\beta_n)_{n\geq1}$ is of stationary type. For
$q\in\N$, let us set $\Omega_q :=
\dis \bigcap_{i\in\cJ}\big( \tau_{Nq+i}=T^*_i\big)\dis \bigcap\big( \beta_q=\beta\big)$. Then, it
is easily seen that $P(\Omega_q)\uparrow 1$ as $q\rw\infty$. For any
$\o \in \Omega_q$, $\theta_{Nq+N+i}(\o)=R^*_i(\o)$, $I^*_{i}(\o)=I_{Nq+N+i}(\o)$ and $\beta_q(\o)=\beta(\o)$. Next let $i\in \cJ$ be fixed.
\begin{equation}\label{eqlemme41}\begin{array}{l} \dis{ \lim_{n\rw \infty} }{\bf E}[
 U^{Nn+N+i}_{\beta_n\wedge\theta_{Nn+N+i}}]= 
 \dis{ \lim_{n\rw \infty} }{\bf E}[
 U^{Nn+N+i}_{\beta_n\wedge\theta_{Nn+N+i}}\{1_{\Omega_n}+1_{\Omega_n^c}\}]  =\dis{ \lim_{n\rw \infty}}{\bf E}[
 U^{Nn+N+i}_{\beta_n\wedge\theta_{Nn+N+i}}
 1_{\Omega_n}]  \\ \\ 
 ={\bf E}[X^{i,\{i\}}_\beta 1_{\{\beta <R_i^*\}}+(X^{i,I_i^*\cup \{i\}}_{R_i^*}\vee 
X^{i,I_i^*}_{R_i^*})1_{\{\beta \ge R_i^*, R_i*<\infty\}}+X^{i,{\cal J}}_\infty 
1_{\{\beta =R_i^*=\infty\}}].\end{array}
\end{equation}
For any $j\in \cJ\setminus{\{i\}}$, we set $\sigma_j=T_j^*$, $\sigma_i=\beta$ and $R=R_i^*\wedge \beta$. By definition of $J_i$, it holds that: \\ 
$$\begin{array}{ll}
{\bf
E}[J_i(T_1^*,...,T_{i-1}^*,\beta,T_{i+1}^*,...,T_N^*)]\\ \\
=   \dis{ \sum_{ I\in\cP }} {\bf E}[\Big\{X^{i,I}_{R}1_{\cap_{j\in
I}\{\sigma_j=R \}\cap \cap_{j\in I^c}\{\sigma_j>R\}}
\Big\}1_{\{R<\infty\}}]+ {\bf E}[
X^{i,\cJ}_{\infty}1_{\{R=\infty\}}]
\\ \\
 = {\bf E}[ X^{i,\{i\}}_{\beta}1_{\{\beta<R_i^* \}}+
X^{i,\cJ}_{\infty}1_{\{R=\infty\}} ]\\ \\\qquad \qquad
+ \dis{ \sum_{I\in\cP, i\notin I} } {\bf E}[\Big\{X^{i,I}_{R_i^*}1_{\cap_{j\in
I}\{\sigma_j=R_i^* \}\cap \cap_{j\in I^c\backslash\{i\}}\{\sigma_j>R_i^*\}}
\Big\}1_{\{R_i^*<\beta\}}]\\ \\\qquad \qquad +
\dis{ \sum_{I\in\cP, i\notin
I} } {\bf E}[\Big\{X^{i,I\cup\{i\}}_{R_i^*}1_{\cap_{j\in I}\{\sigma_j=R_i^*
\}\cap \cap_{j\in I^c\backslash\{i\}}\{\sigma_j>R_i^*\}}
\Big\}1_{\{R_i^*=\beta<\infty\}}]  \\

\\ = {\bf E}[ X^{i,\{i\}}_{\beta}1_{\{\beta<R_i^*\}}+
X^{i,\cJ}_{\infty}1_{\{R=\infty\}}] \\\\
\qquad +  \dis{ \sum_{I\in\cP,
i\notin I}  } {\bf E}[\Big\{X^{i,I}_{R_i^*}1_{\cap_{j\in
I}\{\sigma_j=R_i^*\}\cap \cap_{j\in
I^c\backslash\{i\}}\{\sigma_j>R_i^*\}}
\Big\}1_{\{R_i^*\leq\beta\}}1_{\{R_i^*<\infty\}}]
\\ \\\qquad +  \dis{ \sum_{I\in\cP, i\notin
I} } {\bf E}[\Big\{\Big(X^{i,I\cup\{i\}}_{R_i^*}-X^{i,I}_{R_i^*}
\Big)1_{\cap_{j\in I}\{\sigma_j=R_i^* \}\cap \cap_{j\in
I^c\backslash\{i\}}\{\sigma_j>R_i^*\}}
\Big\}1_{\{R_i^*=\beta<\infty\}}].
\end{array}$$
The last equality is obtained by using 
\begin{equation} \label{decomp}
 X^{i,I}_{R_i^*} 1_{\{R_i^*\le \beta\}\cap \{R_i^*<\infty\}}= 
 X^{i,I}_{R_i^*} 1_{\{R_i^* < \beta\}}+ 
 X^{i,I}_{R_i^*}
1_{\{R_i^*= \beta<\infty\}},
\end{equation}
with the last term in (\ref{decomp}) which is added in the second sum taken over all $I\in\cP$ such that $i\notin I $ and substracted in the last term.
Next for any $I\in\cP$ such that $i\notin I $, it holds
 $$ \{I_{i }^*=I\} = \cap_{j\in I}\{\sigma_j=R_i^*\}\cap
\cap_{j\in
I^c\backslash\{i\}}\{\sigma_j>R_i^*\}.
$$
Therefore 
$$\begin{array}{l}
{\bf
E}[J_i(T_1^*,...,T_{i-1}^*,\beta,T_{i+1}^*,...,T_N^*)]\\\\={\bf E}[ X^{i,\{i\}}_{\beta}1_{\{\beta<R_i^*\}}+X^{i,\cJ}_{\infty}1_{\{R=\infty\}}+\\\\\qq \qq X^{i,I_i^*}_{R_i^*}
1_{\{R_i^*\leq\beta\}}1_{\{R_i^*<\infty\}}+(X^{i,I_i^*\cup\{i\}}_{R_i^*}-X^{i,I_i^*}_{R_i^*}
)1_{\{R_i^*=\beta<\infty\}}].\\\\
  \end{array} $$
Next, substracting $\dis{ \lim_{n\rw \infty} {\bf
E}[ U^{Nn+N+i}_{\beta_n\wedge\theta_{Nn+N+i}}]}$ from the previous equality and using relation (\ref{eqlemme41}) one obtains
  $$ \begin{array}{l}
{\bf
E}[J_i(T_1^*,...,T_{i-1}^*,\beta,T_{i+1}^*,...,T_N^*)]\\\\ 
 = \dis{ \lim_{n\rw \infty}} {\bf
E}[ U^{Nn+N+i}_{\beta_n\wedge\theta_{Nn+N+i}}]\\\\\qq -
{\bf
E} [(X^{i,I^*_{i}}_{R^*_{i}}-X^{i,I^*_{i}\cup\{i\}}_{R^*_{i}} )^+ 1_{\{R^*_{i}=\beta<\infty\}}
+(X^{i,I^*_{i}}_{R^*_{i}}-X^{i,I^*_{i}\cup\{i\}}_{R^*_{i}})^- 1_{\{R^*_{i}<\beta\}} ].
\end{array}
$$ Thus the desired equality (\ref{ebauche1}) is proved which ends the proof.\qed
\begin{lem}\label{lem42}For any $i\in \cJ$ and $\theta\in \cT_0$, we have
\be \label{ineqlam43}\begin{array}{ll} {\bf
E}[J_i(T^*_1,T^*_2,\cdots,T^*_{i-1},\theta,T^*_{i+1},\cdots,
T^*_{N})] +\\\\
  \qq {\bf
E} [\bigg(X^{i,I^*_{i}}_{R^*_{i}}-X^{i,I^*_{i}\cup\{i\}}_{R^*_{i}} \bigg)^+ 1_{\{R^*_{i}=\theta<\infty\}}
+\bigg(X^{i,I^*_{i}}_{R^*_{i}}-X^{i,I^*_{i}\cup\{i\}}_{R^*_{i}} \bigg)^- 1_{\{R^*_{i}<\theta\}}] \\
\\ \leq \varepsilon+{\bf
E}[J_i(T^*_1,T^*_2,\cdots,T^*_{i-1},T^*_i,T^*_{i+1},\cdots,
T^*_{N})]
+\\\\
 \qq  {\bf
E}[ \bigg(X^{i,I^*_{i}}_{R^*_{i}}-X^{i,I^*_{i}\cup\{i\}}_{R^*_{i}} \bigg)^+ 1_{\{R^*_{i}=T^*_i<\infty\}}
+  \bigg(X^{i,I^*_{i}}_{R^*_{i}}-X^{i,I^*_{i}\cup\{i\}}_{R^*_{i}} \bigg)^- 1_{\{R^*_{i}<T^*_i\}}].
\end{array}\ee
\end{lem}
\noindent {\bf Proof}: Let $i\in \cJ$ and $\theta \in \cT_0$. Since $W^{Nn+N+i}$ is a supermartingale, $W^{Nn+N+i}\ge U^{Nn+N+i}$ and $(W^{Nn+N+i}_{k \wedge \mu_{Nn+N+i} })_{k\geq 0}$ is a martingale then
$$\begin{array}{l}
 \dis{ \lim_{n\rw \infty} } {\bf
E}[ U^{Nn+N+i}_{\theta\wedge\theta_{Nn+N+i}}] \leq  \dis{\lim_{n\rw \infty}}{\bf
E}[ W^{Nn+N+i}_{\theta\wedge\theta_{Nn+N+i}}] \leq \dis{ \lim_{n\rw \infty} }{\bf E}
[W^{Nn+N+i}_0]\\\\\qq = \dis{ \lim_{n\rw \infty}}{\bf
E}[ W^{Nn+N+i}_{\mu_{Nn+N+i}}]\\\\\qq \leq \eps+
\dis{ \lim_{n\rw \infty}}{\bf
E}[ U^{Nn+N+i}_{\mu_{Nn+N+i}}] =\eps+ \dis{ \lim_{n\rw \infty} }{\bf
E}[ U^{Nn+N+i}_{\t_{Nn+N+i}\wedge \theta_{Nn+N+i}}].\end{array}
$$
Note that those limits exist due to the almost stationarity of all the decreasing sequences of stopping times which are involved. Finally by Lemma \ref{lem1}, we obtain the desired result since $\lim_{n\rw \infty}\t_{Nn+N+i}=T_i^*$. \qed
 
 \begin{lem}\label{lem43} For any $i_1,i_2\in \cJ$ such that $ i_1\neq i_2$ 
 $$ P\big( T^*_{i_1}=T^*_{i_2}=R^* <\infty\big) = 0. $$
 
 \end{lem}
  {\bf Proof}: Let $\Omega_q:=
\displaystyle\bigcap_{i\in\cJ}\big( \tau_{Nq+i}=T^*_i\big)$ then
  $$\begin{array}{ll}
  \bp\big( T^*_{i_1}=T^*_{i_2}=R^*<\infty \big)& =\bp\big( T^*_{i_1}=T^*_{i_2}=R^*<\infty;\Omega_q \big)+\bp\big( T^*_{i_1}=T^*_{i_2}=R^*<\infty; \Omega_q^c \big)\\ & \leq \bp\big(\tau_{Nq+i_1}=\theta_{Nq+i_1}<\infty \big)+\bp\big( \Omega_q^c \big)\\ & \leq  \bp\big( \Omega_q^c \big).
  \end{array}$$
The second inequality stems from Proposition \ref{prop32}. Taking now the limit w.r.t $q$ completes the proof since  $\lim_{q\rw \infty}\bp\big( \Omega_q^c \big)=0$. \qed

\begin{lem} \label{lem45}For any $i\in \cJ$  $$\begin{array}{ll}
  \{R^*_{i}<T^*_i \}   = \dis{ \bigcup_{j\neq i} }\{T^*_j=R^* <\infty\}.
\end{array}$$
and for all $j\neq i$ we have on $\big(T^*_j=R^*<\infty\big)$
$$ I^*_i=\{j\}. $$
Therefore, under Assumption $(A)$, we have 
\be \label{vanishingterm}
{\bf
E}[\Big(X^{i,I^*_{i}}_{R^*_{i}}-X^{i,I^*_{i}\cup\{i\}}_{R^*_{i}}
\Big)^-1_{\{R^*_{i}<T^*_i \}} ]=0.
\ee
\end{lem}
\ms
{\bf Proof}: First note that  
$$\begin{array}{ll}
 \{R^*_{i}<T^*_i\} =  \dis{ \bigcup_{I\in\cP,
i\notin I} } \bigg(\cap_{j\in I}\{T^*_j=R^*<\infty \}
\cap\cap_{j\in I^c}\{R^*<T^*_j \}\bigg).
\end{array}$$
But by Lemma \ref{lem43} for $I\in \cP$ such that $i\notin I$ and $|I|>1$ we have 
$$\bp\Big(\cap_{j\in I}\{T^*_j=R^* <\infty\}\Big)
  =0,$$
as there cannot exist two different indices $i_1$ and $i_2$ such that $T_{i_1}^*=T_{i_2}^*=R^*<\infty$. Therefore 
$$
\{R^*_{i}<T^*_i\} =  \bigcup_{j\neq i} \{T^*_j=R^* <\infty\}.
$$
Next let $j\neq i$. On the set $\{T_j^*=R^*<\infty\}$, $j\in I_i^*$. Besides, if there exists $i_1\in I_i^*$ satisfying $i_1\neq j$ then one would have $T_j^*=T_{i_1}^*=R^*<\infty$. But this latter is of probability $0$. Thus such an $i_1$ does not exist and $I_i^*=\{j\}$. Finally 
$$ \begin{array}{l}{\bf
E}[\Big(X^{i,I^*_{i}}_{R^*_{i}}-X^{i,I^*_{i}\cup\{i\}}_{R^*_{i}}\Big)^-1_{\{R^*_{i}<T^*_i \}} ]= \dis{  \sum_{j\neq i} } {\bf
E}[\Big(X^{i,I^*_{i}}_{R^*_{i}}-X^{i,I^*_{i}\cup\{i\}}_{R^*_{i}}
\Big)^-1_{\{T^*_j=R^*<\infty \}} ]\\\\
=\dis{ \sum_{j\neq i} } {\bf
E}[\Big(X^{i,\{j\}}_{R^*_{i}}-X^{i,\{i,j\}}_{R^*_{i}}
\Big)^-1_{\{T^*_j=R^*<\infty \}} ]=0
.\end{array}$$
The proof is now complete. \qed
\ms
\begin{rem} \label{rem45} As a by product of Lemma \ref{lem43} we first obtain
 $$\forall i\in \cJ,\,\, \bp\big( T^*_{i}=R^*_{i}<\infty \big) = 0. $$
 Combining this with (\ref{vanishingterm}) in Lemma \ref{lem45}, we deduce that the two last terms in the right-hand side of inequality (\ref{ineqlam43}) in Lemma \ref{lem42} are equal to zero.
 Note that Assumption (A) is crucial to justify (\ref{vanishingterm}).
  
\end{rem}

As a by-product of Lemma \ref{lem42}, Lemma \ref{lem45} and Remark \ref{rem45}, we obtain
the main result of this paper. 
\begin{theorem} \label{mainthm}
The  $N$-tuples $(T_i^*)_{i=1,\dots,N}$ is an $\eps$-{\bf NEP} for the nonzero-sum Dynkin game associated with the payoffs $(J_i)_{i=1,\dots,N}$ of (\ref{J12}).  \qed
\end{theorem}

\indent As a particular case, we end this section by considering a non-zero sum game with $N$ players in discrete time but with finite time horizon $T$ which could be
 random as well.\\ For clarity, we introduce
some extra notations. We fix $T$ in $\N^* $ and, for each $i$ in $\cJ := \{1, \cdots, N\}$ and $I $ in $\cP$, we introduce a collection $( \tilde{X}_t^{i,I})_{t =0,\cdots, T}$  of payoff processes associated with player $i$. We suppose that ($ \tilde{X}_t^{i,I}$)$_{t =0,\cdots, T}$ satisfies 
\begin{equation}\label{newAssumptA}
\displaystyle{ \forall \; t \in \{0,\cdots,T\},\;\; \forall \; (i,j) \in \{1, \cdots, N\}, i \neq j , \quad    
 \tilde{X}_t^{i,\{i,j\}} \le  \tilde{X}_t^{i,\{j\}},
}
\end{equation}
 which is again and by abuse referred as Assumption (A). We also assume  
$$\displaystyle{ 
\forall \; i \in \cJ, \; \forall \; I \in \cP, \quad \quad \tilde{X}_{T}^{i,I} = \tilde{X}_T^{i, \cJ}. } $$
 In such a finite horizon setting, 
this is a common assumption which means that, if the game ends at time $T$, the coalition necessarily consists of all players.    
Next and as in (\ref{J12}), the reward functional $\tilde{J}_i$ for Player $i$ associated with a given $N$-tuple ($T_1, \cdots T_N$) of stopping times valued in $\{0, \cdots, T\}$  is given by
\be \label{newpayoff}  \displaystyle{
\tilde{J}_i(T_1, \cdots T_N)  := \sum_{I\in\cP} \tilde{X}^{i,I}_{R} \mathbf{1}_{ \cap_{ j\in
I}\{T_j = R \}\cap \cap_{j\in I^c} \{T_j> R \}}, 
\;\;\bp -a.s.}
\ee
In what follows, we denote by $\cT_T^N$
(resp. $\cT_T$) the set of all $N$-tuples of stopping times valued in $\{0,\cdots ,T\}$  (resp. the set of all stopping times $\tau$ valued in $\{0,\cdots ,T\}$). 
The following conventions are assumed:\\
(i) the stopping time $R$ satisfies  $R = \mbox{min} \{T_1, \cdots, T_N\}$ and $R$ belongs 
to $\cT_T$;\\
(ii) in analogy with the case with infinite horizon, we impose that if $R(\omega) =T$, then necessarily $I =\cJ$ (or equivalently, the coalition of players consists of all players if the game is stopped at terminal time $T$).\\
Assertion (ii) is satisfied since, by definition of $R$ and on the set $\{\omega,\; R(\omega) =T\}$ one has  $T_i(\omega ) = T $, for any $i$. Then 

 \be \label{equality_timeT} \displaystyle{\tilde{J_i}(T_1, \cdots T_n)\mathbf{1}_{R = T}  = \tilde{J}_i(T, \cdots, T) \mathbf{1}_{R = T}  = X_T^{i, \cJ}}\mathbf{1}_{R = T}.  \ee
\noindent Setting $\overline{\N} := \N \cup \{+\infty\}$, we introduce a  collection  $(X_t^{i,I})_{ t \in \overline{\N}}$ of payoff processes in order to embed the finite horizon case in the infinite horizon one. More precisely, 

\be  \ba{lll} \label{newreward}
\displaystyle{  \forall \; t \in \overline{\N}, \; \forall \;I \subset \mathcal{P},} & \displaystyle{X_t^{i,I} } & \;\displaystyle{ = \;
\tilde{X}_{t }^{i,I},  \; \; \; \textrm{if} \; 0 \le t \le T-1 }; \\  
  \;\;   &\quad \quad \;\;   & \; 
 \displaystyle{= \; \tilde{X}_T^{i, \cJ},  \; \; \; \textrm{if} \;  t\ge T }. \\
 \ea
 \ee
 Let comment the second equality in (\ref{newreward}). Fixing $i \in \cJ$ and $t$ in $\overline{\N} $ such that $t \ge T$, the process $X_t^{i,I}$ does not depend any more on the coalition $I$. 
Finally, we introduce the reward processes $J_i$ associated with player $i$. For any $(T_1, \cdots, T_N) $ in $\cT^N$,
 
 \be \label{auxil_payoff}  \displaystyle{
J_i(T_1, \cdots T_N)  := \sum_{I\in\cP} X^{i,I}_{R} \mathbf{1}_{ \cap_{ j\in
I}\{T_j=R \}\cap \cap_{j\in I^c} \{T_j> R \}},\,\bp -a.s.}, 
\ee
with the stopping time $R$ such that $R = \mbox{min} \{T_1, \cdots, T_N\}$.
We claim that: 

\begin{cor} \label{cor46}
$\\$
 The nonzero-sum game (with $N$ players) with infinite time horizon and reward processes $(J_i)_{i \in \cJ } $ introduced in (\ref{auxil_payoff}) satisfies:
 \begin{enumerate}
  \item[(i)] the collection ($X_t^{i,I}$)$_{ t \in \overline{\N}}$ introduced in (\ref{newreward}) satisfies Assumption ($\bf{A}$);
  \item[(ii)] the same procedure as described in Section 3 provides:
   \begin{itemize}
      \item[(a)] $N$ non-increasing sequences ($\tau_{Nq+i}$)$_{q \in \N}$ initialized by $\tau_1= \cdots \tau_N = \infty $ ;
   \item[(b)]  for any $i\in \cJ$, let us set $T_i^* = \lim_q \searrow \tau_{Nq+i}$. Then the $N$-tuple $(T_i^*)_{i \in \cJ}$  of stopping times is an $\varepsilon$-\textbf{NEP} of the game (with reward processes $(J_i)_{i\in \cJ}$).
   \end{itemize}
   \item[(iii)]  The following relationship holds:
 \begin{equation}\label{equal_rewards} \forall \;(T_1, \cdots T_N) \in \cT^N, \;\;\forall \; i \in \{1,\cdots, N \}, \quad  J^i(T_1, \cdots T_N) =  \tilde{J}^i(T_1\wedge T, \cdots, T_N\wedge T ).
 \end{equation}
 Thus, if for any $i\in\cJ$ we set $\tilde{T_i}^* = T_i \wedge T $, then $(T_i^*)_{\in\cJ}$ is in $\cT_T^N$ and it is an $\varepsilon$-\textbf{NEP} of the nonzero sum game with reward $(\tilde{J}_i)_{i \in \{1,\cdots, N\}}$.
 \end{enumerate}

\end{cor}

\noindent For completeness, we check below all the claims in Corollary \ref{cor46}. In view of  
(\ref{newAssumptA}), the first claim (i) is true 
 and thus, the second claim (ii) results from Theorem \ref{mainthm}.
  To prove the equality in (\ref{equal_rewards}), let us fix a $N$-tuple ($T_1, \cdots T_N$) in $\cT^N$. For this, we need to distinguish the following two cases:\\
  (a) If $R \wedge T = \mbox{min}\{T_1 \wedge T, \cdots  T_N \wedge T\} \le T-1$ then, combining the first equality in (\ref{newreward}) and the definitions of $\tilde{J}_i$ (resp. $J_i$) in (\ref{newpayoff}) (resp. in (\ref{auxil_payoff})),  
  it provides the desired equality.\\ 
 (b) If $R \wedge T=T$ (or equivalenty $R \ge T$) then necessarily for all $i$, $T_i \wedge T=T$ and thus, the desired equality results from (\ref{equality_timeT}).\\
\noindent Relying on Claim (ii)(b), on (\ref{equality_timeT}), (\ref{equal_rewards}) and on the definition in (\ref{equilibrium}) of an $\varepsilon$-NEP then, $(\tilde{T}_i^*)_{i \in \{1, \cdots, N\}}$ provides an $\varepsilon$-{\bf NEP} of the game with reward processes ($\tilde{J}_i$)$_{i \in \{1, \cdots, N\}}$ which ends the proof of Corollary \ref{cor46}.
\qed

\section{Illustration of the constructive algorithm}
$\\ $ In this section and through two explicit examples, 
 we describe our constructive algorithm in discrete time and with finite time horizon $T$ in 
$\N^*$ and we illustrate some properties of the obtained $\varepsilon$-Nash equilibria.
  For sake of clarity, we denote by $(T_i^{*})_{i=1}^N$ any given $N$-tuple produced by the algorithm. Relying on Theorem \ref{mainthm} which is our main result, such a $N$-tuple is an $\varepsilon$-{\bf NEP} of the $N$-player game. In addition, the following property (referred as Claim $(\mathbf{C})$ later) 
  holds:
 \[ \begin{array}{l}
 (\mathbf{C})  \quad \quad   \textrm{the} \; \varepsilon{- \bf NEP} \; (T_i^{*})_{i=1}^N  \; \textrm{may depend on the order of the player in the algorithm.}
  \end{array}
\]
 We stress the fact that, in the algorithm, the "order" of each players is fixed at the beginning and each of the $N$ players successively chooses their optimal stopping time.  As a result, the time horizon of the optimal stopping problem depends on the choice of the $N-1$ other players. Thus, the optimal decision of one player may change depending on his/her order in the construction.  


\subsection{First illustrating example}\label{sec_firstexple}

$\\$ \noindent We study a deterministic case with $N=3$
 players and time horizon $T = 2$. For this, we define the deterministic reward processes $(X_t^{i , I})$ for all $i$ in $\cJ  =\{ 1, 2, 3\}$ and any coalition $I$ in $\cP$. In such a case, we have 
$$ \cP = \{ \{1\} , \{2\}, \{3\}, \{1,2\},\{1,3\},\{2,3\}, \{1,2,3\}  \}.$$
 Since $T=2$, we have to fix all rewards for all players at the three dates $n=0,1,2$.\\
 \noindent At time $n=2$, we impose:
 $$\displaystyle{\forall i \in \cJ,  \quad  X_{2}^{i, \cJ} = 0. }
 $$
\noindent At time $n=0$, 
$$\displaystyle{ \forall \; I\in \cP, \quad X_0^{1, I} = X_0^{2,I} = X_0^{3, I}  =\frac{1}{8}.   } $$
With those conventions, Assumption ($\mathbf{A}$) is satisfied at time $n=0$ and $n= 2$.\\

\noindent Next and for clarity, we collect below in a table all payoff ($X_1^{i, I}$) at time $n = 1$.
\\
\vspace{0,2cm}
\\
\begin{tabular*}{0.75\textwidth}{@{\extracolsep{\fill}} | c | c | c | c |c |c|c|c |}
   \hline
   $i\,\downarrow\,\; 
	  / \;I \,\rightarrow$ & $\{1\}$ & $ \{ 2\}$& $\{ 3\}$&$\{ 1,2\}$&$\{ 1,3\}$&$\{ 2,3\}$&$\{ 1,2,3\} $\\\hline
  1&$\frac{1}{2}$&$\frac{1}{4}$&$\frac{1}{2}$&$\frac{1}{4}$&$\frac{1}{2}$&   $\frac{1}{4}$&$\frac{1}{4}$\\
   \hline
     2&$\frac{1}{2}$&$\frac{3}{2}$&$\frac{1}{2}$&$\frac{1}{4}$&$\frac{1}{2}$&   $\frac{1}{4}$&$\frac{1}{2}$\\
   \hline
   3&$\frac{1}{2}$&$\frac{1}{4}$&$\frac{1}{2}$&$\frac{1}{4}$&$\frac{1}{2}$&   $\frac{1}{4}$&$\frac{1}{4}$\\
   \hline
 \end{tabular*}
\vspace{0,2cm}
\\
\noindent It remains to check Assumption (A) at time $n=1$. First and for Player 1, the following conditions are satisfied: $$ \displaystyle{    \frac{1}{4} = X_1^{1, \{1,2\}} \le X_1^{1, \{2\}} = \frac{1}{4}, \;\;  \textrm{and}  \;\;
  \frac{1}{2}=  X_1^{1, \{1,3\}} \le X_1^{1, \{3\}}=\frac{1}{2} }.$$ 
The payoffs of player 1 and 3 being identical (see the first and third lines above) we obtain the same inequalities as above for player 3. Concerning the second player, one has $$  \displaystyle{  \frac{1}{4} =X_1^{2, \{1,2\}} \le X_1^{2, \{1\}} = \frac{1}{2}, \;\;\textrm{and}  \;\;
  \frac{1}{4} =X_1^{2, \{2,3\}} \le X_1^{2, \{3\}} =\frac{1}{2}.
	}  $$

\noindent For clarity, we provide the main steps of our constructive algorithm which we shall use several times below. 
Recall that $\tau_1 = \tau_2 = \tau_3 =T$. Then, whatever $m \ge  4$ such that $m = 3q_m + i_m$, with $i_m$ in $\{1,2,3\}$, the stopping time $\tau_m$ associated with Player $i_m$ satisfies: 
 \be \label{ost_n} \tau_m :=  \mu_m \mathbf{1}_{ \mu_m < \theta_m } + \tau_{m-N}\mathbf{1}_{ \mu_m \ge \theta_m}, \;\;\ee
where both $\theta_m$ and $\mu_m$ are defined as follows:\\
 (i) $\theta_m := \tau_{m-1}\wedge \cdots  \wedge\tau_{m-(N-1)} $ (in particular $\theta_m = \tau_{m-1} \wedge \tau_{m-2}$, if $N=3$).\\
(ii) Introducing the process $U^m$ as follows: 
\begin{equation}\label{interm_rewardproc}
U_s^m = X_s^{i_m, \{i_m\}} \mathbf{1}_{ s <\theta_m  } +    \big(X_{\theta_m}^{i_m, I_m} \vee X_{\theta_m}^{i_m, I_m \cup \{i_m\}}\big)\mathbf{1}_{ \theta_m \le s < T }  +X_T^{i,\cJ}\mathbf{1}_{\theta_m =T}; \end{equation}
(iii) the $\varepsilon$-optimal stopping time $\mu_n$ satisfies
 $$ \mu_m = \inf \{s \ge 0, \; \textrm{s.t} \; W_s^m \le U_s^m + \varepsilon  \}, \; \; \textrm{where} $$ 
\begin{itemize}

\item[(a)] 
 $W^m$ stands for the Snell envelope process associated with $U^m$;  
\item[(b)]  $I_m$ stands for the coalition of players whose labels are in $\cJ \setminus{ \{i_m\}}$ and which make the decision to stop at time $\theta_m$.\\
 
\end{itemize}

 \subsection{The algorithm applied to the example}
 $\\$
 \indent To begin with, let us provide below two Nash equilibria such that the coalition consists of strictly more than one player. We mention that those Nash equilibria cannot be reached by our explicit algorithm. More precisely, we provide below two 0-{\bf NEP}\footnote{By definition, any $0{-\bf NEP}$ is a fortiori a $\varepsilon$ ${-\bf NEP}$. The other way around is not true in general. } associated with
the game introduced above in Section \ref{sec_firstexple}.\begin{itemize}
\item[(a)] The $3$-tuple
 ($T_1, T_2, T_3$) = ($1,2,1$) is a 0-{\bf NEP}. In this case, both the two players 1 and 3 stop the game at time $t=1$ and thus the optimal coalition is $I^{*} = \{1,3\}$.
\item[(b)] The 3-tuple ($T_1, T_2, T_3$) = ($1,1,1$) is another 0-\textbf{NEP}
 with all players choosing to stop at time $t=1$ and thus 
 the associate coalition is $I^*= \cJ$.
                                                                                                                                                                                       \end{itemize}
  
 \indent  Let prove that these two Nash equilibria cannot be reached as soon as we initialize the algorithm by setting $\tau_i =T=2$ for $i= 1, 2, 3$ reminding here that $\tau_i$ is the stopping time associated with Player $i$.\\
\indent To this end, let construct recursively the sequence $(\tau_m)_{m \ge 4}$ and prove that the algorithm provides the $ 0$-{\bf NEP} 
 $(T_1^*, T_2^*, T_3^*) =(1, 2, 2)$. 
  Thanks to (\ref{ost_n}), it holds $$\tau_4 =  \mu_4 \mathbf{1}_{\mu_4 < 2} + \tau_1 \mathbf{1}_{\mu_4 \ge 2}, $$ 
 since, in that case: $\theta_4 = \tau_3 \wedge \tau_2= 2$. 
 By definition of the reward process $U^4$ in (\ref{interm_rewardproc}) which is associated with player 1 (since $ i_4=1$),
 one obtains $U_0^{4} =X_0^{1, \{1\}} =\frac{1}{8}$, whereas $U_1^4 = X_1^{1, \{1\}} = \frac{1}{2}$ and $U_2^4 =X_2^{1, \cJ}=0$. The Snell envelope process $W^4= SN(U^4)$ being a deterministic process (as it is for $U^4$), it satisfies: $$ W_0^4= W_1^4 = \frac{1}{2} \;\; \textrm{and} \;\; W_2^{4} =0.$$
 Since $U_0^4 < W_0^4$ and $W_1^4 = U_1^4,$ the optimal stopping time is $\mu_4 = 1 < 2$ and therefore $\tau_4 = \mu_4 = 1$ and $\theta_5 =\tau_4 \wedge \tau_3 = \tau_4=1$. 
  Similarly and by definition, $\tau_5$ satisfies 
  \begin{equation}\label{tau5}
  \tau_5 =  \mu_5 \mathbf{1}_{\mu_5 < \theta_5 } + \tau_2 \mathbf{1}_{\mu_5 \ge \theta_5}= \mu_5 \mathbf{1}_{\mu_5 < 1} + \tau_2 \mathbf{1}_{\mu_5 \ge 1}.  
  \end{equation}
  By definition of $U^5$ associated with Player 2 and defined in (\ref{interm_rewardproc}), it holds
  $$ \displaystyle{ U_0^{5} =X_0^{2, \{2\} } = \frac{1}{8} \; \textrm{and} \; U_1^5 =X_1^{2,\{1\}} \vee X_1^{2, \{1,2\}} = X_1^{2,\{1\}}=\frac{1}{2}.}$$
  Since $U_0^{5} < U_1^{5}$, it is not optimal to stop before $\theta_5$ which yields $\mu_5=\theta_5=1$. 
    Using (\ref{tau5}), one obtains $ \tau_5=\tau_2=2$ 
    and $\theta_6 = \tau_5 \wedge \tau_4 =\tau_4 =1 $.\\
  Next and using both $ U_0^6 = X_0^{3, \{3\}} = \frac{1}{8} \; \textrm{and}\; U_1^6 = X_1^{3, \{1\}}\vee  X_1^{3, \{1,3\}} = \frac{1}{2},  $
  the same argumentation as above gives $\mu_6 = \theta_6 =1$ and thus $\tau_6=\tau_3=2.$
 Finally and since $\theta_7 =\theta_4$ then $U^7 =U^4$. Player 1 faces the same optimal stopping problem as before and thus $\tau_7=\tau_4=1$. To sum up, we have obtained $\tau_5=\tau_2=2$, $\tau_6=\tau_3 =2$ and  $\tau_7=\tau_4=1$. Thus and for any $n$, $n \ge 2$ and any $i$ in $\{1,2,3\}$, the three sequences $(\tau_{3n+i})_{n\ge 1}$ are now stationary. The 0-{\bf NEP} ($1, 2,2$) is reached and the game is stopped at time 1 by Player 1 (the coalition is $I^* = \{1\}$).\\

  \indent We conclude by illustrating our claim $(\mathbf{C})$. For this, let
  suppose that the new "order" is (2, 3, 1), meaning that ($\tau_{3n+1}$) (resp. ($\tau_{3n+2}$)$_n$ and ($\tau_{3n+3}$)$_n$) stands for the sequence of stopping times associated with Player 2 (resp. with Player 3 and Player 1).\\  
  \noindent  Once again, we initialize the algorithm by fixing $\tau_i= T=2$ for $i =1,2,3$ and we identify $\tau_4$, $\tau_5 $ and $\tau_6$ recursively defined by (\ref{ost_n}). We first claim that $\tau_4 = \mu_4 =1$. By definition of $U^4$ and since $\theta_4 =\tau_3 =2$,   
    $$   U_0^4 = X_0^{2, \{2\}} = \frac{1}{8}, \; \; U_1^4= X_1^{2, \{2\}} = \frac{3}{2} \;\; \textrm{and}\;\; U_2^4 = X_2^{2, \cJ} =0,$$
    which implies $W_0^4 =W_1^4 = \frac{3}{2}$ and $W_2^4=0 $, and thus $\mu_4 = 1 =\tau_4$.
    Since $\theta_5 = \tau_4 \wedge \tau_3 =\tau_4=  1$, $U_0^5 = X_0^{3,\{3\}} = \frac{1}{8}$ and
    $ U_1^5 =  X_1^{3,\{2\}} \vee X_1^{3,\{3,2\} } =\frac{1}{4} $, the second player (Player 3) has no interest to stop before $\theta_5=1$ and thus $\mu_5 =\theta_5=1$ and $\tau_5 =\tau_2 = 2.$ 
Since $\theta_6 =\tau_5 \wedge \tau_4 =\tau_4=1$, it yields $U_0^6 =X_0^{1, \{1\}} = \frac{1}{8}$ and $ U_1^{6} = X_1^{1, \{2\}} \vee  X_1^{1, \{1, 2\}} = X_1^{1, \{2\}} = \frac{1}{4}$. As above,   
 $\mu_6 =\theta_6=1$ which yields
$\tau_6 =\tau_3=2$ and $\theta_7 = \tau_6 \wedge \tau_5 =2$.\\
 Since $U^7 =U^4$, Player 2 faces the same optimal stopping problem (with horizon $\theta_7 = 2$) meaning that $ W^7=W^4$. The same argumentation as for $\tau_4$ gives $\tau_7 = \tau_4 =1$.
 Thus, for any $i $ in $1, 2,3$, the three sequences ($\tau_{3q+i}$)$_{q \ge 1} $ are stationary, which provides the {\bf NEP} $(T_1^*, T_2^*, T_3^*) =( 2,1 ,2)$ with coalition $I^* = \{2\}$ (consisting of Player 2). \\
   
   
\subsection{Second example with random payoffs}

 \indent We now consider an example of a nonzero-sum game in discrete time with $N = 2$ players which has random reward processes and finite horizon $T=3$.
 We first introduce a Brownian motion $B =(B_n)_{n \in \N^*}$ and an independent sequence of i.i.d.\footnote{\textit{i.i.d} is the standard abbreviation for independent and identically distributed.} random variables $(N_n)_{n \ge 1}$ with common law the uniform law on $\{ -1, \; 1\}$.\\
Let assume that the horizon time $T$ is deterministic and equal to 3. We  introduce below the (random) reward processes associated with each players. For the first player (referred later as Player 1), we set
$$ \displaystyle{\forall \; n \in \{1,2,3\}, \;\; X_n^{1, \{1\}} = B_n;
\; \; X_n^{1 ,\;\{1,2\}} = B_n +\frac{1}{2}, \; \; \textrm{and}\; \;  X_n^{1 ,\;\{2\}} =  B_n +1,} $$
whereas for the second player (referred as Player 2), we set
$$ \displaystyle{\forall \; n \in \{1,2,3\}, \;\;X_n^{2, \{2\}} = B_n + N_n;\; \; X_n^{2 ,\;\{1,2\}} = B_n + N_n+\frac{1}{2}; \; \; \textrm{and}\;  X_n^{2 ,\;\{1\}} =  B_n + N_n +1.} $$
Since, for any $n$ in $\{1,\;2,\;3\}$ both conditions $  X_n^{1 ,\;\{1,2\}} \le  X_n^{1 ,\;\{2\}}$ and  $ X_n^{2 ,\;\{1,2\}} \le  X_n^{2 ,\;\{1\}}$ hold, Assumption ($\mathbf{A}$) is satisfied. \\
On such a discrete time setting, we introduce the following filtration $ (\mathcal{F}_n)_{n \ge 1}$ 
$$ \forall \;\; n \in \{1, \cdots T\},\quad  \mathcal{F}_n = \sigma\big( B_i, N_i, \; i \in \{1, \cdots, n\} \big). $$  
From the definitions of $(B_n)_n$ and $(N_n)_n$ and using both the independence and/or martingale properties, we deduce: 
\begin{equation}\label{conditional_esp}
\E \big( B_{n+1}|\mathcal{F}_n\big) = \E\big( B_{n+1}| B_n\big) =B_n  \;\textrm{and} \; \E \big( N_{n+1}|\mathcal{F}_n\big) =
\E \big( N_{n+1}\big) =0.
\end{equation}

\noindent To compute the Snell envelope $W:= SN(U)$ of
the process $U$, we recall its (backward recursive) construction in discrete time  
\begin{equation}\label{snellenvelope}
  W_T =U_T \; \; \textrm{and}\;\;  W_n = \mbox{max} \{ U_n;\; \E \big(W_{n+1} |\mathcal{F}_n \big) \},\; \textrm{for} \;  n \;=\; T-1, \cdots, 1.
\end{equation}
 \noindent Let apply our algorithm by providing an explicit (random) $\varepsilon$-{\bf NEP},  
$0 \le \varepsilon <\frac{1}{2}$. 
As in the previous paragraph, we construct both the two sequences $(\tau_m )$ and $(\mu_m )_{m\ge 3}$ with the first one initialized as follows $\tau_1 = \tau_2= 3.$ 
Assuming here that Player 1 begins, he/she chooses first its $\varepsilon$-optimal stopping time $\mu^3 = \mu^{3}(\varepsilon, \omega)$(\footnote{From now, we omit both symbols $\varepsilon$ and $\omega$: contrary to the first example, all stopping times $\theta_m$, $\mu_m$ and $\tau_m$ are a priori random and so it is for the {\bf NEP}.}) defined as follows
\begin{equation}\label{eps_optimtime}
\mu_3  = \mbox{Inf}\{ n \in \{1,2,3\},\; W_n^3 \le U_n^3 +\varepsilon  \}, 
\end{equation}
where, as in (\ref{interm_rewardproc}),  $W^3 = SN(U^3)$ and $U^3$ satisfies
\begin{equation}\label{intermpayoff_P1}
\forall \; n \in \{1,2,3\}, \q 
U_n^3 = \underbrace{X_n^{1, \{1\}} }_{\displaystyle{  =B_n} }\mathbf{1}_{n < 3} + \underbrace{X_3^{1, \{1,2\} }}_{\displaystyle{ = B_3+\frac{1}{2}} } \mathbf{1}_{n=3}, \;\; \textrm{since} \; \theta_3=\tau_2=3. 
\end{equation}
Using both (\ref{snellenvelope}) and the martingale property of $(B_n)$ stated in (\ref{conditional_esp}), we obtain \\ $\forall \; n \in  \{ 1,\;2,\;3\} \;\; W_n^3= B_n+\frac{1}{2}$. Thus and since $\varepsilon < \frac{1}{2}$, the $\varepsilon$-optimal stopping time $\mu_3$ defined in (\ref{eps_optimtime}) satisfies $\mu_3=\theta_3 =3$. 
By definition of $\tau_3$ in (\ref{ost_n}) and since $N=2$, it holds
$$\tau_3 = \mu_3\mathbf{1}_{ \mu_3 < \theta_3} +\tau_1 \mathbf{1}_{\mu_3=\theta_3}= \mu_3\mathbf{1}_{ \mu_3 <  3} + \tau_1 \mathbf{1}_{\mu_3 =3} = \tau_1 = 3, \; 
$$
 which implies $ \theta_4= \tau_3 = 3  $. Next, $\mu_4$ satisfies
$$ \mu_4  = \mbox{Inf}\{ 1\le n \le 3,\; \;
W_n^4 \le U_n^4 +\varepsilon  \}, $$
with the reward process $U^4$ such that
\begin{equation}\label{intermpayoff_P2}
\displaystyle{
 U_n^4 = \underbrace{X_n^{2,\{2\}} }_{= \displaystyle{B_n+N_n}} \mathbf{1}_{n <3 } +\underbrace{X_n^{2, \{1,2\}}}_{ \displaystyle{ = B_{3}+N_{3}+ \frac{1}{2}} } \mathbf{1}_{ n =3
}. }
\end{equation}
Again using (\ref{ost_n}), $\tau_4$ is such that: $\tau_4 = \mu_4 \mathbf{1}_{\mu_4 < \theta_4} + \tau_2   \mathbf{1}_{\mu_4 = \theta_4}
 $.\\
\noindent 
To identify the (random) $\varepsilon$-{\bf NEP}, let compute
the Snell envelope process $W^4 =SN(U^4)$ associated with $U^4$ expressed in (\ref{intermpayoff_P2}).\\
Since $\theta_4 = \tau_3 = 3$ and using (\ref{intermpayoff_P2}), we claim  \begin{equation}\label{U4case2} U_1^4= X_1^{2,\{2\}} =B_1+N_1, \; U_2^4= B_2+N_2, \;\;
U_3^4 = B_3 +N_3 +\frac{1}{2}.
\end{equation} By definition of $W^4$ in (\ref{snellenvelope}), $  W_{3}^4=U_3^4 =B_3 +N_3 +\frac{1}{2}$. Using both (\ref{conditional_esp}) and (\ref{snellenvelope}),
$$ W_2^4 = \mbox{max}\{    
B_2+N_2; \underbrace{\E\big( B_3 |
\mathcal{F}_2 \big)}_{ = \E(B_3 | B_2) =B_2 }+\underbrace{\E \big( N_3 |\mathcal{F}_2 \big)}_{ = \E(N_3)} +\frac{1}{2}  \}=  \mbox{max} \{ B_2+N_2; B_2+\frac{1}{2} \}. $$
This leads to
\begin{equation}\label{snellenv_date2}
 W_2^4  = (B_2+ 1)\mathbf{1}_{N_2=1} + (B_2+ \frac{1}{2})\mathbf{1}_{N_2=-1}.
\end{equation} 
Finally 
$ W_1^4 = \mbox{max} \{  B_1+N_1; \E(W_2^4 |\mathcal{F}_1)  \},$ with
$$ \E(W_2^4 |\mathcal{F}_1)  = (B_1+1) \mathbf{P}(N_2=1) +  (B_1+ \frac{1}{2})\mathbf{P}(N_2=-1)=B_1 + \frac{3}{4}, $$
which yields
\begin{equation}\label{snellenv_date1}
W_1^4  = (B_1+1) \mathbf{1}_{N_1=1} + (B_1 + \frac{3}{4})\mathbf{1}_{N_1=-1}.  
\end{equation}
\noindent
It remains to distinguish the following three cases:\\

\noindent \underline{Case (i):} $N_1  =1:$\\
Relying on (\ref{snellenv_date1}), one has $W_1^4 = (B_1+1) \mathbf{1}_{N_1=1} =U_1^4 $, which gives $\mu_4 =1$. Thus $\tau_4 = \mu_4= 1$ and $ \theta_5 = \tau_4 = 1$.\\ 
Next since condition $\mu_5 < \theta_5=1$ cannot hold, and since $\tau_5 = \mu_5 \mathbf{1}_{\mu_5 < \theta_5} + \tau_3\mathbf{1}_{\mu_5 \ge \theta_5} $
then, necessarily $\tau_5=\tau_3 = 3 $ and $\theta_6=\tau_5= 3$. It now suffices to prove
 that $\tau_6=\tau_4 \;( =1)$. Since $\theta_6 =\theta_4=3$, Player 2 again solves the same optimal stopping problem with reward process $U^6$ equal to $U^4$. Thus $W^6 =W^4 $ which leads to $\mu_6 = \mu_4 =1$ and implies $\tau_6 = \mu_6 =\mu_4=\tau_4$. The desired claim $\tau_6=\tau_4$ is established. As a result, both sequences $(\tau_{2q+1})_{q \ge 1}$ and $(\tau_{2q+2})_{q \ge 1}$ are now stationary. Thus on $\{N_1 =1\}$, the $\varepsilon$-{\bf NEP} $(T_1^*, T_2^*) = (3,1)$ is reached.\\
\vspace{0,2cm}
\\ 
\noindent \underline{Case (ii):} $N_1=-1$ and $N_2 =1$:\\
In this case and in view of (\ref{snellenv_date2}) and (\ref{snellenv_date1}), we claim 
\begin{equation} \label{snellenv-caseii}
 W_1^4 =B_1 + \frac{3}{4} > B_1-1+\varepsilon =U_1^4 + \varepsilon \; \; \textrm{and} \;\; W_2^4 = B_2+1=U_2^4.
\end{equation}
 This gives $\mu_4 = 2$ and $\tau_4=2 = \theta_5$. 
Since $ \theta_5= \tau_4 = 2$, Player 1 faces an optimal stopping problem with $\varepsilon$-stopping time $\mu_5$ and with reward process $U^5$ given by: 
$$U_n^5 =X_n^{1,\{1\}}\mathbf{1}_{n < \theta^5} +X_{\theta_5}^{1,\{2\}} \vee X_{\theta_5}^{1, \{1,2\}}\mathbf{1}_{\theta_5 \le n <2}  = B_n \mathbf{1}_{n <2} + (B_2+1)\mathbf{1}_{n = 2}.$$
From the martingale property of $B_n$ and since $\varepsilon < \frac{1}{2}$, it implies $\mu_5=\theta_5 =2$ and $\tau_5 = \tau_3 = 3$. 
It remains to prove that $ \tau_6 =\tau_4 =2$ so that, as in case (i) above, both $(\tau_{2q+1})_{q \ge 1}$ and $(\tau_{2q+2})$ are stationary. 
Since $\theta_6 = 3 = \theta_4$ then $U^6 =U^4$ with $U^4$ given in (\ref{U4case2}) which yields $W^4 =W^6$. 
Therefore, $ \mu_6 = \mu_4= 2$ and $\tau_6=\mu_6 = 2 =\tau_4$, which is the desired claim.  \\
 On $\{N_1 =-1;\;N_2= 1 \}$, we obtain the $\varepsilon$-{\bf NEP} $(T_1^*, T_2^*) = (3, 2)$ .\\
\vspace{0,2cm}
\\
\underline{Case (iii):} $N_1=-1$ and $N_2 =-1$:\\
On this last case and since $\varepsilon < \frac{1}{2}$, it holds 
$$ W_1^4 =B_1 + \frac{3}{4} > B_1-1+\varepsilon =U_1^4 + \varepsilon \q \; \textrm{and} \q\; W_2^4= B_2+\frac{1}{2}   > B_2+ N_2 +\varepsilon = U_2^4 + \varepsilon ,$$ 
which means that $\mu_4 = 3$ and $\tau_4= 3$. We thus obtain $\theta_5 =\tau_4 =3$. Since $\theta_3=\theta_5 =3$, then $U^5=U^3$ and thus $W^5= W^3$. Player 1 faces the same optimal stopping problem as before, which yields $\tau_5 =\tau_3=3$. Once again and for any $i=1,2$, $(\tau_{2q+i})_{q \ge 1}$ are stationary sequences and we obtain the $ \varepsilon$-{\bf NEP} $(T_1^*, T_2^*) = (3,3)$.\\

\noindent We provide in a final remark two last comments concerning the constructive algorithm.

\paragraph*{Remark} \begin{itemize}
 \item[(1)] Let fix $\varepsilon$ such that $\varepsilon \ge \frac{1}{2}$ and let suppose that the constructive algorithm begins with Player 1. We already know that $W_n^3 =B_n +\frac{1}{2}$ for any $n \in \{1,2,3\}$. Therefore, $W_1^3 = B_1+\frac{1}{2} \le B_1+ \varepsilon =  U_1^3+ \varepsilon $ and the $\varepsilon$-optimal stopping time is $\mu_3=1$. Using (\ref{ost_n}), we obtain $\tau_3  =\mu_3 =1$, which implies $\theta_4= \tau_3 =1$. Since $\mu_4 <\theta_4$ cannot hold, we obtain $\tau_4=\tau_2=3$ and $\theta_5= \tau_4=3$. 
Noting that $\theta_5=\theta_3$ as in
 Case (iii) above, we obtain $U^5=U^3$ which implies that $W^5=W^3$ and $\tau_5 = \mu_5 = \mu_3 =1$.   
 Thus, for any $ i =1,2$, $(\tau_{2q+i})_{q\ge 1}$ are stationary and the $ \varepsilon $-{\bf NEP} 
$(T_1^*, T_2^*) = (1,3)$ is reached.

\item[(2)] On the contrary and when $\varepsilon <\frac{1}{2} $, 
 we illustrate that the order between players in the constructive algorithm does not change the $ \varepsilon $-{\bf NEP}.\\
If the first player in the algorithm is Player 2 then its reward process $\bar U^3$ 
 is defined similarly as $U^4$ in (\ref{U4case2}) and he faces a stopping problem with  horizon $\theta_3 =\tau_1= 3$ .
 Therefore, Player 2 stops either at time
 $ \mu_3 = 1$, $2$ or $3$ as in case (i)-(iii) above. Next, Player 1 faces an optimal stopping problem with horizon
  $\theta_4=\mu_3 $ and reward process $\bar U^4$ defined as follows
 $$ \displaystyle{ \bar U_n^4 =  B_n \mathbf{1}_{n < \theta_4} + 
 ( B_n + \frac{1}{2} )\mathbf{1}_{ n \ge \theta_4 }. } $$ 
Thus, whatever $\theta_4$ and using once again the martingale property of $B$, Player 1 has never interest to stop before $\theta_4$ which implies that 
$ \mu_4 = \theta_4 =1$. 
The obtained $ \varepsilon $-{\bf NEP} is the same as in case (i)-(iii) above, which proves the desired claim.
\end{itemize}

To sum up, we have highlighted through two explicit examples some properties of the constructive algorithm. In particular we show that it can produce several $\eps$-Nash equilibria. More precisely we illustrate that the order between the $N$ players in the construction may influence the reached $ \eps$-{\bf NEP}.
 We also note that, due to Assumption $(\bf{ A})$ and for any $\eps$-{\bf NEP} produced by the constructive algorithm, the coalition $I^*$ of players which decide to stop the game consists of at most one player when the game terminates before the horizon time. As a result, any other $\eps$-{\bf NEP} such that $|I^*| > 1$  
 cannot be obtained by using this algorithm.

\end{document}